\documentclass[12pt]{amsart}

\usepackage{graphicx}
\usepackage{amsthm}
\swapnumbers

\theoremstyle{plain}
\newtheorem{Thm}{Theorem}

\newtheorem{Coro}[Thm]{Corollary}
\newtheorem{Lem}[Thm]{Lemma}

\theoremstyle{definition}
\newtheorem{Def}[Thm]{Definition}

\begin{document}

\title{On tunnel number one knots \\ that are not $(1,n)$}
\author{Jesse Johnson and Abigail Thompson}
\address{\hskip-\parindent
        Department of Mathematics\\
        Yale University\\
        New Haven, CT 06511\\
        USA}
\email{jessee.johnson@yale.edu}
\address{\hskip-\parindent
        Department of Mathematics\\
        University of California\\
        Davis, CA 95616\\
        USA}
\email{thompson@math.ucdavis.edu}
\subjclass{Primary 57M}
\keywords{Tunnel number one, bridge number, curve complex}

\thanks{This research was supported by NSF VIGRE grant 0135345}

\begin{abstract}
We show that the bridge number of a $t$ bridge knot in $S^3$ with respect to an unknotted genus $t$ surface is bounded below by a function of the distance of the Heegaard splitting induced by the $t$ bridges.  It follows that for any natural number $n$, there is a tunnel number one knot in $S^3$ that is not $(1,n)$.
\end{abstract}

\maketitle

\section{Introduction}

A compact, connected, closed, orientable surface $S$ embedded in 
$S^3$ is \textit{standardly embedded} if the closure of each 
component of its complement is a handlebody.  Equivalently, $S$ is a 
Heegaard surface for $S^3$.  A knot $K$ is in $n$-bridge position 
with respect to $S$ if the intersection of $K$ with each handlebody 
is a collection of $n$ boundary parallel arcs.

For $n \geq 1$, we will say that $K$ is $(t,n)$ if $K$ can be put in 
$n$-bridge position with respect to a standardly embedded, genus $t$ 
surface $S$.  We will say that $K$ is $(t,0)$ if $K$ can be isotoped 
into $S$.  If $K$ is $(t,n)$ for some $n$ then $K$ is $(t,m)$ for 
every $m \geq n$.  Thus the important number is the smallest $n$ such 
that $K$ is $(t,n)$.

A set of arcs properly embedded in the the complement of a knot $K$ 
is an \textit{unknotting system} if the complement of a regular neighborhood of
$K$ and the arcs is a handlebody.  The \textit{tunnel number} of 
$K$ is the minimum number of arcs in an unknotting system for $K$. 

Let $K$ be a knot in $S^3$ and $\Sigma$ the Heegaard splitting of the knot complement induced by a $t$-tunnel decomposition for $K$.  Hempel defined a distance $d(\Sigma)$ for Heegaard splittings using the curve complex.  We will prove the following:

\begin{Thm}
\label{boundthm}
If $K$ is $(t,n)$ then $K$ is $(t,0)$ or $d(\Sigma) \leq 2n + 2t$.
\end{Thm}

Every tunnel number $t$ knot is $(t+1,0)$.  The question is for what values 
of $n$ can a tunnel number $t$ 
knot be $(t,n)$.  Moriah and Rubinstein~\cite{mr:neg} showed that 
there exist tunnel number one knots that are $(1,2)$, but not $(1,1)$.  Morimoto, 
Sakuma and Yokota~\cite{MSY} and Eudave-Mu\~noz~\cite{eud:tor} 
constructed further examples of knots that are not $(1,1)$.  Eudave-Mu\~noz has recently announced the existence of tunnel number one knots that are not 
$(1,2)$.  The first author of this paper~\cite{me:bridges} showed that for tunnel number one knots, $d(\Sigma)$ can be arbitrarily large.  Thus Theorem~\ref{boundthm} implies the following:

\begin{Coro}
\label{maincoro}
For every $n \in \mathbf{N}$, there is a tunnel number one knot $K$ such that $K$ is not $(1,n)$.
\end{Coro}

The proof in \cite{me:bridges} is non-constructive and therefore does not provide actual examples of knots with high toroidal bridge number.  Since this note first appeared as a preprint, Minsky, Moriah and Schleimer~\cite{mms:dist} have given a constructive proof that there are $t$-tunnel knots in $S^3$ with arbitrarily high distance splittings.  They conclude, using Theorem~\ref{boundthm}, that for every $t$ and $k$, there is a $t$ tunnel knot that is not $(t,k)$.

We describe weakly incompressible surfaces in Section~\ref{wrsect} and the curve complex in Section~\ref{ccsect}.  Theorem~\ref{boundthm} and Corollary~\ref{maincoro} are proved in Section~\ref{bndsect}.

\section{Weakly Compressible Surfaces}
\label{wrsect}

A properly embedded, two sided surface $S$ in a 3-manifold $M$ is 
\textit{compressible} if there is a disk $D$ in $M$ such that 
$\partial D$ is an essential simple closed curve in $S$ and the 
interior of $D$ is disjoint from $S$.  If $S$ is not compressible 
then $S$ is \textit{incompressible}.

Assume that $S$ separates $M$ into components $X$ and $Y$.  Then $S$ 
is \textit{strongly compressible} if there are disks $D_1$ and $D_2$ 
such that $\partial D_1$ and $\partial D_2$ are disjoint, essential 
simple closed curves in $S$, the interior of $D_1$ is contained in 
$X$ (disjoint from $S$) and the interior of $D_2$ is contained in 
$Y$.  If $S$ is not strongly compressible then $S$ is \textit{weakly 
incompressible}.

A properly embedded surface $S$ is \textit{boundary compressible} if 
there is a disk $D \subset M$ such that $\partial D$ consists of an 
essential arc in $S$ and an arc in $\partial M$.  A separating 
surface $S$ is \textit{strongly boundary compressible} if there are 
boundary compressing disks on opposite sides of $S$ with disjoint 
boundaries, or a boundary compressing disk and a compressing disk on 
opposite sides of $S$ with disjoint boundaries.  A surface is 
\textit{weakly boundary incompressible} if $S$ is not strongly 
boundary compressible and $S$ is not strongly compressible.

\begin{Lem}
\label{wilem}
Let $M$ be a compact 3-manifold and $F$ a closed, separating, 
incompressible torus embedded in $M$. Let $A$, $B$ be the closures of 
the components of the complement of $F$.  Let $S$ be a second surface 
which separates $M$.  If $S \cap A$ is weakly boundary incompressible 
in $A$ and $S \cap B$ is empty or incompressible and boundary 
incompressible in $B$  then $S$ is weakly incompressible in $M$.  If 
$S \cap A$ and $S \cap B$ are both incompressible and boundary 
incompressible,  then $S$ is incompressible in $M$.
\end{Lem}

\begin{proof}
Assume for contradiction $S$ is strongly compressible.  Then there 
are disks $D_1$, $D_2$ properly embedded on opposite sides of $S$ 
such that $\partial D_1 \cap \partial D_2 $ is empty.

Assume $D_1$ and $D_2$ have been chosen transverse to $F$ and with a 
minimal number of  components in $(D_1 \cup D_2) \cap F$.  If $D_1$ 
and $D_2$ are disjoint from $F$ then both disks must be in $A$ 
because $S \cap B$ is incompressible.  This contradicts the 
assumption that $S \cap A$ is weakly boundary incompressible. 
Without loss of generality, assume $F \cap D_1$ is not empty.

Because $F$ is incompressible and any loop in $D_1$ is trivial in 
$D_1$, any loop component of $D_1 \cap F$ must be trivial in $F$. 
Compressing $D_1$ along an innermost such loop will reduce the number 
of components of intersection without changing its boundary.  Thus 
minimality implies $D_1 \cap F$ is a collection of arcs.  Similarly, 
if $D_2 \cap F$ is not empty then $D_2 \cap F$ is a collection of 
arcs.

An outermost arc $\beta$ in $D_1$ cuts off a disk whose boundary 
consists of an arc $\alpha$ in $F$ and an arc $\beta$ in $S \cap A$ 
or $S \cap B$.   If the arc $\beta$ is trivial in $S \cap B$ or $S 
\cap A$ then it can be pushed across $F$ (taking any other arcs with 
it) and reducing $(D_1 \cup D_2) \cap F$.  Thus we can assume that 
$\beta$ is essential in $S \cap A$ or $S \cap B$.

If $\beta$ is in $S \cap B$ then the outermost disk is a boundary 
compression disk for $S \cap B$.  Because $S \cap B$ is boundary 
incompressible, this is not possible so $\beta$ must be in $S \cap A$ 
and $D_1$ contains a boundary compression disk $D$ for $S \cap A$.

If $D_2$ is disjoint from $F$ then $D_2$ is a compression disk for $S 
\cap A$.  This compression disk is on the opposite side from $D$ and 
$\partial D$ is disjoint from $\partial D_2$.  This contradicts the 
assumption that $S \cap A$ is weakly boundary incompressible.  If 
$D_2$ intersects $F$ then, as with $D_1$, an outermost disk argument 
implies that $D_2$ contains a boundary compressing disk $D'$ for $S 
\cap A$.  The disks $D$ and $D'$ are disjoint and on opposite sides 
of $S \cap A$, again contradicting weak boundary incompressibility.

The case in which $S \cap A$ and $S \cap B$ are both incompressible 
and boundary incompressible proceeds similarly, but more 
easily.
\end{proof}

To apply Lemma~\ref{wilem} to knots, we need a result regarding thin position for a 
knot in the 3-sphere with respect to a standard genus g Heegaard splitting.  
The result  follows from unpublished work of C. Feist ~\cite{feist:genus}.  His 
Theorem 5.5 implies:

\begin{Lem}
\label{thinlem}
If a knot $K$ is $(t,n)$ and not $(t,0)$ then either (case 1) there 
is a bicompressible, weakly boundary incompressible meridinal genus 
$t$ surface with at most $2n$ boundary components in the complement 
of $K$ or (case 2) there is an incompressible, boundary 
incompressible meridinal surface with genus at most $t$ and at most 
$2n$ boundary components in the complement of $K$.
\end{Lem}

\section{The Curve Complex}
\label{ccsect}

Let $H$ be a 3-manifold with boundary and let $\Sigma$ be a component 
of $\partial H$.

\begin{Def}
The \textit{curve complex} $C(\Sigma)$ is the graph whose vertices are isotopy
classes of simple closed curves in $\Sigma$ and edges connect vertices
corresponding to disjoint curves.
\end{Def}

For more detailed descriptions of the curve complex, 
see~\cite{hemp:cc} and~\cite{mm:hyp}.

\begin{Def}
The \textit{boundary set} $\mathbf{H} \subset C(\Sigma)$ 
corresponding to $H$ is the set of vertices $\{l \in C(\Sigma) : l$ 
bounds a disk in $H\}$.
\end{Def}

Given vertices $l_1, l_2$ in $C(\Sigma)$, the distance $d(l_1, l_2)$ 
is the geodesic distance: the number of edges in the shortest path 
from $l_1$ to $l_2$.  This definition  extends to a definition of 
distances between subsets $X,Y$ of $C(\Sigma)$ by defining $d(X,Y) = 
\min \{d(x,y) : x \in X, y \in Y\}$ and for distances between a point 
and a set similarly.

Given a compact, connected, orientable 3-manifold $M$ and a compact, 
connected, closed, separating surface $\Sigma$, let $A$ and $B$ be 
the closures of the complement in $M$ of $\Sigma$.  Then $\Sigma$ is 
a component of $\partial A$ and a component of $\partial B$.  Let 
$X$,$Y$ be the boundary sets in $C(\Sigma)$ of $A$ and $B$, 
respectively.  If $X$ and $Y$ are non-empty, we will define $d(\Sigma) = d(X,Y)$.

This situation arises in a knot complement as follows:  Let $M$ be 
the complement of a regular neighborhood of a knot $K$ in $S^3$ and 
let $\tau_1,\dots,\tau_t$ be a collection of properly embedded arcs 
in $M$.  The arcs $\tau_1,\dots,\tau_t$ are called a collection of 
\textit{unknotting tunnels} for $K$ if the complement in $M$ of a 
regular neighborhood $N$ of $\bigcup \tau_i \cup \partial M$ is a 
handlebody.  Let $\Sigma$ be the boundary component of the closure of 
$N$ that is disjoint from $\partial M$.   The surface $\Sigma$ 
separates $M$ and allows us to define $d(\Sigma)$ as above.  For $t = 
1$, Lemma 4 and Lemma 11 of \cite{me:bridges} imply the following 
Lemma:

\begin{Lem}
\label{highdistlem}
For every $N$, there is a knot $K$ in $S^3$ and an unknotting tunnel 
$\tau$ such that for $\Sigma$ constructed as above $d(\Sigma) > N$.
\end{Lem}

In \cite{me:bridges}, it is shown that $d(\Sigma)$ bounds below both 
the bridge number of $K$ and the Seifert genus of $K$. 
Theorem~\ref{boundthm} provides a similar bound for the toroidal bridge 
number.

\section{Bounding Distance}
\label{bndsect}

A compact, separating surface $\Sigma$ properly embedded in a 
manifold $M$ is called \textit{bicompressible} if there are 
compressing disks for $\Sigma$ in both components of $M \setminus 
\Sigma$.

Given a bicompressible, weakly incompressible surface $\Sigma$, let 
$A$, $B$ be the closures of the complements of $M \setminus \Sigma$. 
If we compress $\Sigma$ into $A$, the resulting surface, $\Sigma'$, 
separates $A$.  It may be possible to compress $\Sigma'$ still 
further into the component of $A \setminus \Sigma'$ which does not 
contain $\Sigma$, creating a new surface which again separates $A$.

Let $\Sigma_A$ be the result of compressing $\Sigma'$ away from 
$\Sigma$ repeatedly, until the resulting surface has no compression 
disks on the side which does not contain $\Sigma$.  Let $\Sigma_B$ be 
the result of the same operation, but compressing $\Sigma$ maximally 
into $B$.   Define $\Sigma^*$ to be the submanifold of $M$ bounded by 
$\Sigma_A$ and $\Sigma_B$.  Following \cite{tom:dist} (with slightly 
different notation), we will say that weakly incompressible surfaces 
$\Sigma$ and $S$ are \textit{well separated} if $S^*$ can be isotoped 
disjoint from $\Sigma^*$.  We will say that $\Sigma$ and $S$ are 
parallel if $S$ can be isotoped to be parallel to $\Sigma$.  The 
following is Theorem 3.3 in \cite{tom:dist}.

\begin{Thm}[Scharlemann and Tomova~\cite{tom:dist}]
\label{tomthm}
If $\Sigma$ and $S$ are bicompressible, weakly incompressible, 
connected, closed surfaces in $M$ then either $\Sigma$ and $S$ are 
well separated, $\Sigma$ and $S$ are parallel, or $d(\Sigma) \leq 2 - 
\chi(S)$.
\end{Thm}

This theorem is the key to the following proof.  Note that $2 - \chi(S)$ is precisely twice the genus of $S$.

\begin{proof}[Proof of Theorem~\ref{boundthm}]
Let $M$ be the complement in $S^3$ of a neighborhood of a knot $K$ 
and assume $K$ is $(t,n)$.  By Lemma~\ref{thinlem}, there is either 
an incompressible, boundary incompressible or a bicompressible, weakly boundary 
incompressible $2k$-punctured genus $t$ surface $T$ properly embedded 
in $M$ with $k \leq n$.  

Let $M'$ be the complement in $S^3$ of a 
neighborhood of the connect sum of $k$ trefoil knots.

There is a collection $T'$ of $k$ pairwise disjoint, properly 
embedded, essential annuli in $M'$ and there is a homeomorphism $\phi 
: \partial M \rightarrow \partial M'$ which sends $\partial T$ onto 
$\partial T'$.  Let $M''$ be the result of gluing $M$ and $M'$ via 
the map $\phi$.  The image in $M''$ of $T' \cup T$ is a closed, genus 
$t+k$ surface which we will call $S$.  The Euler 
characteristic of $S$ is $2 - 2(k + t)$.  Because $T$ is incompressible 
or weakly incompressible and $T'$ is incompressible, 
Lemma~\ref{wilem} implies that $S$ is either incompressible or weakly 
incompressible.

Lemma~\ref{wilem} also implies that the image in $M''$ of $\Sigma$ is 
weakly incompressible because $\Sigma$ is weakly incompressible in 
$M$ and $\Sigma \cap M'$ is empty.  

Suppose  $T' \cup T$ is 
compressible but weakly incompressible.  Then  by 
Theorem~\ref{tomthm}, either $\Sigma$ and $S$ are parallel, the 
surfaces are well-separated or $d(\Sigma) \leq 2(k+t) \leq 2n + 2t$. 
To complete the proof of this case we will show that $\Sigma$ and $S$ 
are not parallel or well separated.

First we will show that the surfaces are not parallel.  The surface 
$\Sigma$ bounds a submanifold containing the closed, incompressible 
torus $\partial M$.  If $\Sigma$ and $S$ are parallel then the 
complement of $S$ contains an incompressible torus $A$, isotopic to 
$\partial M$.  Assume for contradiction this is the case.  Any loop 
in the intersection $A \cap \partial M$ must be trivial in both 
surfaces or essential in both, as both surfaces are incompressible. 
Any trivial loop of intersection can be eliminated by an isotopy of 
$A$ which keeps $A$ disjoint from $S$, so we can assume $A \cap S$ is 
empty or consists of essential loops.

If $A \cap S$ is empty then $A$ is contained in $M$ or $M'$.  If $M$ 
contains an essential torus then as noted in \cite{thm:dcp}, 
$d(\Sigma) \leq 2$ and we are done.  Thus we will assume the only 
incompressible surface in $M$ is boundary parallel.  Such a surface 
cannot be disjoint from $T \subset S$.

Each component of the complement in $M'$ of $T'$ is homeomorphic to 
an unknot complement or a trefoil knot complement.  Thus an 
incompressible surface in $M'$ which does not intersect $T'$ bounds 
an unknot complement or a trefoil complement. If $\partial M$ is 
isotopic to one of these surfaces, then $M$ must be an unknot or 
trefoil complement.  In either case, $d(\Sigma) \leq 2$ 
(see~\cite{me:bridges}).  Thus we will assume $A \cap S$ must be 
non-empty.

Let $A'$ be a component of $A \cap M$.  An incompressible annulus 
properly embedded in $M$ is always boundary parallel, so one 
component of $M \setminus A'$ is a solid torus.  The surface $S$ 
cannot be contained in this solid torus, so $A'$ can be isotoped 
across $\partial M$, reducing $A \cap \partial M$.  This implies $A$ 
is disjoint from $\partial M$, which we saw above is a contradiction. 
Hence $A$ and $\Sigma$ are not parallel.

To show that the surfaces are not well separated, consider the 
subsets $\Sigma^*$ and $S^*$ of $M''$  defined above.  The surface 
$\Sigma$ compresses down to a ball on one side and to a neighborhood 
of $\partial M$ on the other, so we can take $\Sigma^*$ to be the 
image in $M''$ of $M$.  If $\Sigma$ and $S$ are well separated  then 
$S$ can be isotoped out of $M''$. After the isotopy, $\partial M''$ 
is an incompressible surface in the complement of $S$.  Thus there is 
an incompressible torus, isotopic to $\partial M$ in the complement 
of $S$.  We showed that no such surface exists, so $\Sigma$ and $S$ 
are not well separated.  

Now suppose  $T' \cup T$ is incompressible. 
The arguments of Theorem~\ref{tomthm} apply to this case as well, 
although considerably simplified by the fact that  $T' \cup T$ is 
incompressible instead of weakly incompressible.  The details of this 
case are left to the reader.    
 \end{proof}

\begin{proof}[Proof of Corollary~\ref{maincoro}]
By Lemma~\ref{highdistlem}, there is a knot $K$ with unknotting 
tunnel $\tau$ such that for the induced Heegaard splitting $\Sigma$, 
$d(\Sigma) > 2n + 2$.  As noted in~\cite{me:bridges}, every 
unknotting tunnel for a torus knot has distance at most 2, so $K$ is 
not $(1,0)$.  Thus by Theorem~\ref{boundthm}, $K$ is not $(1,n)$.
\end{proof}

\bibliographystyle{abbrv}
\bibliography{bridges2}

\end{document}